\documentclass[12pt]{amsart}

\usepackage{amsmath,amssymb}
\usepackage{amsfonts}
\usepackage{graphicx}
\usepackage{epsfig}
\input{amssym.def}

\theoremstyle{plain}
\newtheorem{proposition}{Proposition}[section]
\newtheorem{lemma}[proposition]{Lemma}

\newtheorem{definition}[proposition]{Definition}

\newtheorem*{corollaire }{Corollaire}
\newtheorem{theorem}{Theorem}

\newtheorem{remark}[proposition]{Remark}

\newcommand{\diam}{\mbox{\rm diam}}

\def\C{\mathbb C}

\def\P{\mathbb P}
\def\Q{\mathbb Q}

\def\Z{\mathbb Z}

\def\cB{{\mathcal B}}
\def\cD{{\mathcal D}}

\def\cO{{\mathcal O}}

\def\hm{\widehat{m}}
\def\hr{\widehat{r}}

\def\trho{\widetilde{\rho}}

\def\fm{{\frak m}}

\def\proof{\par\noindent {\it Proof.} }
\def\lemma{\par\noindent {\bf Lemma.} }
\def\proposition{\par\noindent {\bf Proposition.} }
\def\remark{\par\noindent {\bf Remark.} }

\def\definition{\par\noindent {\bf Definition.} }

\def\pp{{{\mathbb P}(\C_p)}}

\begin{document}
\title{Wild recurrent critical points.}

\author{Juan Rivera-Letelier}
\address{Departamento de Matem\'aticas, Universidad Cat\'olica del Norte, Casilla 1280 - Antofagasta, Chile.}
\date{February 2004, revised June 2004.}
\email{rivera-letelier@ucn.cl}
\begin{abstract}
It is conjectured that a rational map whose coefficients are algebraic over $\Q_p$ has no wandering components of the Fatou set.
R.~Benedetto has shown that any counter example to this conjecture must have a wild recurrent critical point.
We provide here the first examples of rational maps whose coefficients are algebraic over $\Q_p$ and that have a (wild) recurrent critical point.
In fact, we show that there is such a rational map in every one parameter family of rational maps that is defined over a finite extension of $\Q_p$ and that has a Misiurewicz bifurcation.
\end{abstract}
\maketitle

\section{Introduction.}
Fix a prime number $p$, let $\Q_p$ be the field of $p$-adic numbers and let $\C_p$ be the completion of an algebraic closure of $\Q_p$.
We consider rational maps whose coefficients in $\C_p$ are algebraic over $\Q_p$, viewed as a dynamical systems acting on the projective line $\pp$.

\subsection{No wandering domains conjecture.}
Fix a rational map $R \in \C_p(z)$ whose coefficients are algebraic over $\Q_p$.
Just as in the complex case, the projective line $\pp$ is partitioned into two sets: the {\it Fatou set} where the dynamics is stable and the {\it Julia set} where the dynamics is (topologically) expanding~\cite{Hs}.
The Fatou set is open and dense in $\pp$ and it is in turn partitioned into three open sets, where the dynamics of $R$ is: {\it contracting}, (eventually) {\it quasi-periodic} and {\it wandering}~\cite{these}.

The situation here is analogous to the complex case, but Sullivan has shown that complex rational maps do not have wandering domains~\cite{Su}.
By analogy to the complex case, it is conjectured that the wandering part of the Fatou set is empty~\cite{Benowandering},~\cite{these}.
In its simplest formulation this conjecture is as follows.

\

{\par\noindent {\bf Conjecture.} {\it Every wandering disk is attracted to an attracting cycle.}}

\

\subsection{Wild recurrent critical points.}
R.~Benedetto has shown that any counter example to the no wandering domains conjecture must have a non-periodic wild recurrent critical point~\cite{Benowandering}.
A {\it critical point} of a rational map $R$ is a point in $\pp$ where the local degree of $R$ is strictly bigger than~1.
If moreover this local degree is a multiple of $p$, then the critical point is called {\it wild}.
A point of $\pp$ is {\it recurrent} under $R$ if it is accumulated by its forward orbit under $R$.

We provide here the first examples of rational maps {\it whose coefficients are algebraic over $\Q_p$}, having a (wild) recurrent critical point that is not periodic.
Such a critical point must belong to the Julia set, so for such rational maps there is a delicate balance between the local contraction of the critical point and the topological expansion of the Julia set.

\

{\par\noindent {\bf Theorem~A.}}
{\it In every one parameter family of rational maps that is defined over a finite extension of $\Q_p$ and that has a Misiurewicz bifurcation, there is a rational map whose coefficients are algebraic over $\Q_p$ and that has a non-periodic recurrent critical point.}

\

A one parameter family of rational maps $R$ has a {\it Misiurewicz bifurcation} at a parameter $t_0$, if there is a critical point that is mapped to a repelling periodic point by some iterate $R^\ell_{t_0}$ of $R_{t_0}$, but for $t \neq t_0$ near $t_0$ this does not happen.
See Section~\ref{Misiurewicz bifurcations} for precise definitions.

We give an algebraic condition that guarantee that the rational map in the theorem has coefficients in a given finite extension of $\Q_p$, see Corollary~2 in Section~\ref{Misiurewicz cascade}.
For a given $d > 1$, we show that the family,
$$
P_t(z) = - p^{-d}z (z - 1)^d + t,
\mbox{ for } t \in \C_p
$$
has a Misiurewicz bifurcation at $t = 0$ and that it satisfies this algebraic condition.
We conclude that there are arbitrarily small parameters $t$ in $\Q_p$, for which the critical point~1 of $P_t$ is recurrent but not periodic.
Note that the local degree of $P_t$ at~1 is equal to $d$, so when $p$ divides $d$ this critical point is wild.
\subsection{Comparison with the complex case.}
In the complex case C.~McMullen has shown that every non-trivial {\it algebraic} family of rational maps has bifurcations, except for some special families of Latt\`es maps~\cite{Mc1}.
Moreover Misiurewicz bifurcations are dense in the bifurcation locus~\cite{Mc2}.
So Misiurewicz bifurcations occur in virtually every algebraic family of complex rational maps.

This statement does not hold in the ultrametric setting considered here.
To begin with, the family of quadratic polynomials $P_c(z) = z^2 + c$, for $c \in \C_p$, is a non-trivial algebraic family of non-Latt\`es rational maps, and yet it does not have Misiurewicz bifurcations.
In fact, when $|4c| \le 1$ all cycles are non-repelling and when $|4c| > 1$ the (unique) finite critical point escapes to infinity, so it cannot be mapped to a repelling periodic point.
Besides parabolic bifurcations, that have a small impact in the overall dynamics, the quadratic family seems to be {\it locally} stable.

\

{\par\noindent {\bf Problem~1.}}
{\it Characterize those algebraic families having a Misiurewicz bifurcation.}

\

A {\it marked critical point} of a one parameter family of rational maps $R$, is a function $c : t \mapsto c_t$, such that for every parameter $t$ the point $c_t$ is a critical point of $R_t$.
A marked critical point is {\it active at a parameter~$t_0$} if on any neighborhood of $t_0$ the family $\{ R^n(c) \}_{n \ge 1}$ is not equicontinuous.
A natural candidate for the {\it bifurcation locus} is the set of those parameters for which there is a (locally) active critical point, compare with Proposition~2.4 of~\cite{Mc2}.

\

{\par\noindent {\bf Problem~2.}}
{\it Consider a one parameter family of rational maps $R$ with a marked critical point that is active at some parameter $t_0$.
Is $t_0$ accumulated by Misiurewicz bifurcations?}

\subsection{Remarks.}
As $\pp$ is totally disconnected (every connected component is reduced to a point), is not natural to consider the decomposition of the Fatou set into its connected components, as it is done in the complex setting.
Moreover the Fatou set is usually equal to $\pp$!
There is however a natural (dynamical) partition of the Fatou set, but it is not easy to describe.
This is done in detail in~\cite{Fatou}, see also~\cite{components} and~\cite{these}.
The situation is much more natural when one considers the action of the rational map on the Berkovich analytic space of $\pp$.

Recently Benedetto has shown examples of (polynomial) rational maps whose coefficients are {\it transcendental over $\Q_p$}, having a wandering domain~\cite{Bewandering}.
Using similar techniques, it can be shown that near every Misiurewicz bifurcation of a wild critical point there are (transcendental) parameters having a wandering domain.

Here we follow in part the axiomatic treatment of Misiurewicz bifurcations in~\cite{Mc2}.
\subsection*{Acknowledgments.}
I'm grateful to J.~Kiwi and R.~Benedetto for useful comments on a preliminary version of this article.
This paper was written during a visit of the author at IMPA and it was partially supported by MECESUP~UCN0202 and PROSUL.

\section{Preliminaries.}
Fix a prime number $p$, let $\Q_p$ be the field of $p$-adic numbers and let $\C_p$ be the completion of an algebraic closure of $\Q_p$.
{\it All finite extensions of $\Q_p$ that we will consider will be subfields of $\C_p$}.
We denote by $| \cdot |$ the norm of $\C_p$, by $\C_p^*$ its multiplicative group and
$$
|\C_p^*|
=
\{ |z| \mid z \in \C_p^* \}
=
\{ p^r > 0 \mid r \in \Q \}.
$$ 
We also denote $\cO_p = \{ z \in \C_p \mid |z| \le 1 \}$ the ring of integers of $\C_p$ and $\fm_p = \{ z \in \C_p \mid |z| < 1 \}$ the maximal ideal of $\cO_p$.

For $a \in \C_p$ and  $r \in |\C_p^*|$ we call
$$
\{ |z - a| \le r \}
\mbox{ or }
\{ |z - a| < r \}
$$
{\it open} or {\it closed ball of $\C_p$}, respectively.
When $r > 0$ does not belong to $|\C_p^*|$ these sets coincide and they form what we call an {\it irrational ball of $\C_p$}.
By definition a ball of $\C_p$ is irrational if and only if $\diam(B) \not \in |\C_p^*|$, so {\it the diameter of a open or closed ball of $\C_p$ belongs to $|\C_p^*|$}.

We identify the projective line $\pp$ of $\C_p$ with $\C_p \cup \{ \infty \}$.
So, for every subfield $K$ of $\C_p$ the corresponding projective line $\P (K) \subset \pp$ is identified with $K \cup \{ \infty \}$.
The {\it chordal metric} $\Delta$ on $\pp$, is defined in homogeneous coordinates by
$$
\Delta([x, y], [x', y'])
=
\frac{|xy' - x'y|}{\max \{ |x|, |y| \} \cdot \max \{ |x'|, |y'| \}}.
$$
When $z, z' \in \C_p$ we have
$$
\Delta(z, z')
=
\frac{|z - z'|}{\max \{ 1, |z| \} \cdot \max \{ 1, |z'| \} } \le |z - z'|.
$$
\section{Holomorphic functions.}\label{holomorphic functions}
Let $n$ be a positive integer.
A {\it closed} or {\it open polydisk} in $\C_p^n$ is a product $n$ of closed or open balls of $\C_p$, respectively.
When $n = 2$, polydisks are also called {\it bidisks}.

A {\it holomorphic function} defined on a polydisk $\underline{D}$ in $\C_p^n$, is a function given by a convergent power series of the form
$$
(z_1, \ldots, z_n)
\mapsto
\sum_{\nu_1, \ldots, \nu_n \ge 0} c_{\nu_1, \ldots, \nu_n} \cdot (z_1 - \zeta_1)^{\nu_1} \cdot \ldots \cdot (z_n - \zeta_n)^{\nu_n}.
~~\footnote{
Usually holomorphic functions are defined over more general domains.
As our results are of local nature, holomorphic functions defined over polydisks will be enough for our purposes.}
$$
Such a function has a power series expansion at each point $(\zeta_1, \ldots, \zeta_n)$ in $\underline{D}$, that converges in all of $\underline{D}$.
The sum, the product and the composition of holomorphic functions are holomorphic functions.

For a finite extension $K$ of $\Q_p$, we say that a holomorphic function $f$ defined on a polydisk $\underline{D}$ is {\it defined over $K$}, if the power series of $f$ at some point of $\underline{D}$ in $K^n$, has coefficients in $K$.
In that case the power series of $f$ at each point of $\underline{D}$ in $K^n$ has coefficients in $K$.
Clearly a holomorphic function $f$ defined over $K$ is also defined over every finite extension $L$ of $K$.
So the image by $f$ of an element of $L$ belongs to $L$.

The sum, the product and composition of holomorphic functions defined over $K$ are also holomorphic functions defined over $K$.

\subsection{}\label{basic}
The following lemma states some basic properties of holomorphic functions that we will use frequently.
For the proof in the case of the unit polydisk see Section~5.1 of~\cite{BGR}.
The general case follows easily from this case.

\

\lemma
{\it Let $D_1, \ldots, D_n$ be closed balls of $\C_p$ of diameters $r_1, \ldots, r_n$, respectively.
Then, for $(\zeta_1, \ldots, \zeta_n) \in \underline{D} = D_1 \times \ldots \times D_n$, a power series
$$
f(z_1, \ldots, z_n)
=
\sum_{\nu_1, \ldots, \nu_n \ge 0} c_{\nu_1, \ldots, \nu_n} \cdot (z_1 - \zeta_1)^{\nu_1} \cdot \ldots \cdot (z_n - \zeta_n)^{\nu_n}
$$
converges on $\underline{D}$ if and only if
$$
|c_{\nu_1, \ldots, \nu_n}| \cdot r_1^{\nu_1} \cdot \ldots \cdot r_n^{\nu_n} \to 0
\mbox{ as }
\max \{ \nu_1, \ldots, \nu_n \} \to \infty.
$$
In that case
$$
\sup_{\underline{D}} |f|
=
\sup_{\nu_1, \ldots, \nu_n \ge 0} |c_{\nu_1, \ldots, \nu_n}| \cdot r_1^{\nu_1} \cdot \ldots \cdot r_n^{\nu_n}.
$$}
\subsection{Lemma}\label{injective holomorphic functions}(Injective holomorphic functions){\bf .}

{\it Given $r \in |\C_p^*|$ set $D = \{ |z| \le r \}$ and let $f(z) = c_0 + c_1 z + \ldots$ be a power series converging on $D$.
Then the following assertions hold.
\begin{enumerate}
\item[1.]
$f$ is injective on $D$ if and only if for every $k > 1$ we have $|c_k| r^{k - 1} < |c_1|$.
In that case
$$
f(D) = \{ z \in \C_p \mid |z - c_0| \le |c_1| r \}.
$$
\item[2.]
If $f$ is injective, then its inverse function $g : f(D) \to D$ is holomorphic.
\item[3.]
If $f$ is injective and defined over a finite extension $K$ of $\Q_p$, then the inverse of $f$ is defined over $K$ and $f$ induces a bijection between $D \cap K$ and $f(D) \cap K$.
\end{enumerate}}
\proof
Part~1 is an easy consequence of Schwarz' Lemma, as stated in~\cite{these}, Section~1.3.1.
To prove part~2, we reduce to the case $c_0 = 0$ and $c_1 = 1$, so that $f(D) = D$.
Since $f$ converges on $D$ we have $|c_k| r^{k - 1} \to 0$ as $k \to \infty$.
By part~1, the hypothesis that $f$ is injective implies that for $k > 1$ we have $|c_k| r^{k - 1} < 1$.
Therefore $\gamma = \sup \{ |c_k| r^{k - 1} \mid k > 1 \} < 1$.

Let us now define inductively $b_i \in \C_p$, in such a way that
$$
\sum_{1 \le j \le i} b_j(f(z))^j \equiv z \mod z^{i + 1}.
$$
Set $b_1 = 1$, so that this is satisfied for $i = 1$.
Suppose that this holds for some positive integer $i$ and let $a \in \C_p$ be defined by
$$
\sum_{1 \le j \le i} b_j(f(z))^j \equiv z + a z^{i + 1} \mod z^{i + 2}.
$$
Then $b_{i + 1} = -a$ satisfies the inductive condition.

By definition $b_i$ is a linear combination, with coefficients in $\Z$, of terms of the form 
$$
b_j \cdot \prod_{1 \le \ell \le j} c_{k_\ell},
$$
where $j \ge 1$ and $k_1 + \ldots + k_j = i$, with $k_\ell \ge 1$.
So it follows by induction that $b_i \in \Z[c_1, \ldots, c_i]$.
On the other hand,
$$
\left| b_j \cdot \prod_{1 \le \ell \le j} c_{k_\ell} \right| r^i
\le 
|b_j|r^j \cdot \min \{ \gamma^j, |c_{k_1}| r^{k_1 - 1}, \ldots, |c_{k_j}| r^{k_j - 1} \}.
$$
It follows by induction on $i$ that $|b_i| r^i < r$, for $i > 1$ and that $|b_i| r^i \to 0$ as $i \to \infty$.
Therefore the power series $z + b_2 z^2 + \ldots$ converges on $D$ to the inverse of $f$.

For part~3, just observe that $b_i \in \Z[c_1, \ldots, c_i]$, so the inverse of $f$ is defined over $K$ whenever $f$ is.
The last assertion of part~3 follows from the fact that holomorphic functions defined over $K$ map points of $K$ to points of $K$.
$\hfill \square$
\subsection{Inverse branches.}\label{inverse branches}
Let $f : D \to \C_p$ be a holomorphic function.
If for a point $z_0 \in D$ we have $f'(z_0) \neq 0$, then $f$ is injective on a small ball containing $z_0$ and therefore $f$ has a local inverse $g$ at $z_0$ (Lemma~\ref{injective holomorphic functions}).
The function $g$ is holomorphic and it is defined over a finite extension $K$ of $\Q_p$, whenever $f$ is defined over $K$ and $z_0 \in K$.
\subsection{Lemma.}\label{Rouche}
{\it Let $f$ be an injective holomorphic function defined on a closed ball $D$ of $\C_p$, such that the ball $f(D)$ contains $0$.
Then for every holomorphic function $g : D \to \C_p$ satisfying 
$$
\sup_D |g| < \diam(f(D)),
$$
there is $z_0 \in D$ such that $f(z_0) = g(z_0)$.
Moreover, if $f$ and $g$ are defined over a finite extension $K$ of $\Q_p$, then $z_0 \in K$.}
\proof
After an affine coordinate change on the domain we assume $D = \cO_p$.
Since $f(z) = c_0 + c_1 z + \ldots$ is injective, for every $i > 1$ we have $|c_i| < |c_1|$ and $f(\cO_p) = \{ |z - c_0| \le |c_1| \}$ (Lemma~\ref{injective holomorphic functions}).
So the hypothesis $0 \in f(\cO_p)$ implies $|c_0| \le |c_1|$.

If we set $g(z) = c_0' + c_1'z + \ldots $, then
$$
\sup \{ |c_i'| \mid i \ge 0 \} = \sup_{\cO_p} |g| < \diam(f(\cO_p)) = |c_1|.
$$
So Lemma~\ref{injective holomorphic functions} applied to $f - g$ implies that the holomorphic function $f - g$ is injective on $\cO_p$ and $(f - g)(\cO_p) = \{ |z| \le |c_1| \}$.
So there is $z_0 \in \cO_p$ such that $(f - g)(z_0) = 0$.

If $f$ and $g$ are defined over $K$, then $f - g$ is defined over $K$ and $z_0 \in K$ (part~2 of Lemma~\ref{injective holomorphic functions}).
$\hfill \square$
\subsection{One parameter families of holomorphic functions.}
A {\it one parameter family of holomorphic functions}, is a holomorphic function defined on a bidisk.
Such a function $f : B \times D \to \C_p$ can be written in the form
$$
f_t(z) = c_0(t) + c_1(t) (z  - \zeta) + c_2(t) (z - \zeta)^2 + \ldots,
$$
for each $\zeta \in D$, where each $c_i$ is a holomorphic function defined on $B$ and for each $t \in B$ the function $f_t$ is holomorphic.

When $f$ is defined over a finite extension $K$ of $\Q_p$ and $\zeta \in K$, the functions $c_i$ are defined over $K$.
In that case, for every $t \in K$ the function $f_t$ is defined over $K$.

\subsection{Lemma.}\label{one parameter injectivity}
{\it
Let $f$ be a one parameter family of holomorphic functions defined on a bidisk containing $(0, 0)$, of the form
$$
f_t(z) = c_0(t) + c_1(t) z + \ldots
$$
Given $r, r' \in | \C_p^* |$ set $B = \{ t \in \C_p \mid |t| \le r \}$ and $D = \{ z \in \C_p \mid |z| \le r' \}$.
Assume that $r > 0$ and $r' > 0$ are such that for every $t \in B$, we have
$$
|c_1(t)| = |c_1(0)| > 0, \ \
|c_0(t)| \le |c_1(0)|r'
$$
$$
\mbox{ and for $i > 1$ } \ \
|c_i(t)| (r')^{i - 1} < |c_1(0)|.
$$
Then for every $t \in B$ the holomorphic function $f_t$ maps $D$ bijectively onto
$$
D' = \{ z \in \C_p \mid |z| \le |c_1(0)| r' \}.
$$
Moreover, the inverse functions $(f_t|_D)^{-1}$ form a one parameter family of holomorphic functions $g$, defined on $B \times D'$.

If $f$ is defined over a finite extension $K$ of $\Q_p$, then so is the family~$g$.}

\proof
Replacing $f$ by $f - f(0)$ we assume $c_0 \equiv 0$.
The hypothesis that for every $t \in B$ we have $|c_1(t)| = |c_0(0)| = 1$, implies that the function $t \mapsto c_1(t)^{-1}$ is holomorphic, cf.~\cite{FvP}~Lemme~I.2.5.
So replacing $f$ by $c_1^{-1} f$ we assume $c_1 \equiv 1$.

Then the proof is similar to the proof of part~2 of Lemma~\ref{injective holomorphic functions}.
$\hfill \square$

\subsection{}\label{one parameter inverse branches}
The following lemma is an easy consequence of the previous one.

\

\lemma
{\it 
Let $f$ be a one parameter family of holomorphic functions defined on a bidisk $\underline{D}$ and let $(t_0, z_0) \in \underline{D}$ be such that $f_{t_0}'(z_0) \neq 0$.
Then for every small ball $D$ containing $z_0$ there are balls $B$ and $D'$ of $\C_p$, such that $t_0 \in B$ and such that for every $t \in B$ the function $f_t$ induces a bijection between $D$ and $D'$.
Moreover, the inverse functions $(f_t|_D)^{-1}$ form a one parameter family of holomorphic functions $g$ defined on $B \times D'$.

If $f$ is defined over a finite extension $K$ of $\Q_p$ and $t_0, z_0 \in K$, then so is the family $g$.}
\subsection{Lemma.}\label{near critical inverse branches}
{\it
Let $\delta \in |\C_p^*|$ satisfying $\delta \le |p|^2$ and consider a one parameter family of holomorphic functions of the form
$$
f_t(z) = c_d(t) z^d + c_{d + 1}(t) z^{d + 1} + \ldots,
\mbox{ with } c_d(0) \neq 0,
$$
defined for $t$ and $z$ close to $0$.
Then for every $t$ and every $\xi$ close to $0$, the function $f_t$ maps the ball
$$
\{z \in \C_p \mid |z - \xi| \le \delta |\xi| \}
$$
bijectively onto the ball
$$
\{ w \in \C_p \mid |w - c_d(0) \xi^d| \le \delta |d \xi^d| \cdot |c_d(0)| \}.
$$
Moreover, for a fixed $\xi$ the corresponding inverse maps form a one parameter family of holomorphic functions.

If $f$ is defined over a finite extension $K$ of $\Q_p$ and if $\xi \in K$, then this family of inverse maps is defined over $K$.}

\proof
Let $B$ be a ball containing $0$ such that for every $t \in B$ we have $|c_d(t) - c_d(0)| \le \delta |d| \cdot |c_d(0)|$.
Then $|c_d(t)| = |c_d(0)|$ and 
$$
\{ w \in \C_p \mid |w - c_d(t) \xi^d| \le \delta |d \xi^d| \cdot |c_d(0)| \}
$$
$$
=
\{ w \in \C_p \mid |w - c_d(0) \xi^d| \le \delta |d \xi^d| \cdot |c_d(0)| \}.
$$
So the function $t \mapsto c_d(t)^{-1}$ is holomorphic on $B$ and, replacing $f$ by $c_d^{-1} f$, we may assume that $c_d(t) \equiv 1$.

Let $\eta \in \C_p$ be of norm equal to $\delta$ and let $r' > 0$ be small enough so that for $j > d$ we have $\sup_B |c_j(t)| (r')^{j - d} < \delta |pd|$.
Set $D = \{ |z| \le r' \}$ and fix $\xi \in D$.
Consider the one parameter family of holomorphic functions $h$, defined on $B \times \cO_p$ by
$$
h_t(z) = \xi^{-d} f_t(\xi (1 + \eta z))  - 1 = b_0(t) + b_1(t) z + \ldots
$$
Observe that, as $|\eta| = \delta \le |p|^2$, the coefficients of the polynomial $(1 + \eta z)^d - 1 - \eta d z$ have norm at most $\delta |pd|$.
Therefore, for every $t \in B$ we have
$$
|b_1(t) - \eta d| \le \delta |pd|
\ \  \mbox{ and for $j \neq 1$ } \ \
|b_j(t)| \le \delta |pd|.
$$
So $h$ satisfies the hypothesis of Lemma~\ref{one parameter injectivity} with $r' = 1$.
The assertions follow easily from the conclusions of this lemma.
$\hfill \square$
\section{Dynamics of rational maps.}
We consider a rational map $R$ with coefficients in $\C_p$ as a dynamical system acting on $\pp$.
We refer the reader to~\cite{Hs} and~\cite{these} for background on dynamics of rational maps.
For a positive integer $n$ we denote by $R^n$ the $n$-th iterate of $R$ and $R^0$ will denote the identity map.
 
Rational maps are locally holomorphic functions, in the sense that they are given by a convergent power series on a neighborhood of each point of $\pp$, in suitable coordinates.
This follows from the fact every rational map has a power series expansion on a neighborhood of every point of $\C_p$ that is not a pole.
\subsection{Periodic points.}\label{periodic points}
A point $z_0 \in \pp$ is {\it fixed} by a rational map $R \in \C_p(z)$, if $R(z_0) = z_0$.
The point $z_0$ is {\it periodic} by $R$ if for some integer $n \ge 1$ the point $z_0$ is fixed by $R^n$.
When $n$ is the least integer with this property, we say that $n$ is the {\it period} of $z_0$.
In that case the derivative $\lambda = (R^n)'(z_0)$, in a coordinate such that $z_0 \neq \infty$, is invariant under coordinate changes and its is called the {\it multiplier} of $z_0$.
We say that $z_0$ is a {\it repelling, indifferent} or {\it attracting periodic point}, according to $|\lambda| > 1$, $|\lambda| = 1$ or $|\lambda| < 1$, respectively.

\

\lemma
{\it Suppose that $a \in \C_p$ is a repelling fixed point of a rational map $R$ and set $\lambda = R'(a)$.
Denote by $g$ the local inverse of $R$ at $a$.
Then there is a ball $D$ containing $a$, such that $g(D) \subset D$ and such that the sequence of functions $\{ \lambda^i(g^i - a) \}_{i \ge 1}$ converges uniformly on $D$ to a holomorphic function $\zeta$.
Moreover $\zeta$ is locally injective at $a$.}

\

By definition the function $\zeta$ satisfies the functional equation $\zeta \circ g = \lambda^{-1} \zeta$, so that $\zeta$ conjugates $g$ to the map $z \mapsto \lambda^{-1} z$ on a neighborhood of $a$.

\

\proof
After an affine coordinate change assume that $a = 0$ and that $g$ is defined on $\cO_p$.
Then $g$ is of the form
$$
g(z) = \lambda^{-1} z(1 + c_1 z + c_2 z^2 + \ldots ),
$$
with $\lim_{i \to \infty} |c_i| = 0$.
Set $C = \sup_{i \ge 1} |c_i|$ and note that for every $z \in \cO_p$ we have $|\lambda g(z) - z| \le C|z|^2$.
So if we set $D = \{ |z| < C^{-1} \}$, then for every $z \in D$ we have $|g(z)| = |\lambda|^{-1}|z| < |z|$.
So $g(D) \subset D$ and it follows by induction that for every $i \ge 1$ we have $|g^i(z)| = |\lambda|^{-i} |z|$.

Set $\zeta_i = \lambda^i g^i$, so that for every $z \in D$ we have $|\zeta_i(z)| = |z|$ and
$$
|\zeta_{i + 1}(z) - \zeta_i(z)| \le C|\lambda|^{i}|g^i(z)|^2 = (C|z|^2)|\lambda|^{-i}.
\hfill \square
$$  
\subsection{Misiurewicz maps and homoclinic orbits.}\label{Misiurewicz and homoclinic}
A {\it rational map with marked critical point} is a pair $(R, c)$ consisting of a rational map $R$ and a critical point $c$ of $R$.
When $R$ has coefficients in a finite extension $K$ of $\Q_p$ and $c \in K$, we say that $(R, c)$ {\it is defined over} $K$.

A rational map with marked critical point $(R, c)$ is {\it Misiurewicz}, if the critical point $c$ is mapped to a repelling periodic point of $R$ under iteration.
In this case, replacing $R$ by an iterate if necessary, we assume that the critical point $c$ is mapped to a repelling fixed point of $R$ under iteration.
Fix such $(R, c)$ and let $\ell \ge 1$ be the minimal integer for which $a = R^\ell(c)$ is a repelling fixed point.
A {\it bi-infinite orbit} is a sequence $\{ z_i \}_{i \in \Z}$ such that for every $i \in \Z$ we have $R(z_i) = z_{i + 1}$.
A {\it homoclinic orbit} is a bi-infinite orbit $\{ z_i \}_{i \in \Z}$ for which
$$
\lim_{|i| \to \infty} z_i = a.
$$
As $a$ is a repelling fixed point, this implies that for some $m \in \Z$ we have $z_m = a$.
Moreover, if we denote by $g$ the local inverse of $R$ at $a$, then for some $n \in \Z$ the point $z_n$ belongs to the domain of definition of $g$ and for $j \ge 0$ we have $z_{n - j} = g^j(z_n)$.
\subsection{Ramification and critical homoclinic orbits.}\label{critical homoclinic}
A point $z$ is an {\it unramified preimage} of a point $w$, if for some positive integer $n \ge 1$ we have $R^n(z) = w$ and if the local degree of $R^n$ at $z$ is equal to~1.
A point is {\it unramified} if it has infinitely many unramified preimages.

Let $(R, c), a, \ell, \ldots$ be as in \S\ref{Misiurewicz and homoclinic} above.
We say that a homoclinic orbit $\{ z_i \}_{i \in \Z}$ is {\it critical}, if $z_0 = c$ and if for every $j < 0$ the point $z_j$ is an unramified preimage of $c$.
In this case we have $z_j = a$ for every $j \ge \ell$.

\

\proposition
{\it Let $(R, c)$ be a Misiurewicz rational map with marked critical point, defined over a finite extension $K$ of $\Q_p$.
If the critical point $c$ is unramified, then there is a critical homoclinic orbit contained in a finite extension of $K$.}
\proof
Let $\ell \ge 1$ be the least integer such that $R^\ell(c)$ is a repelling periodic point of $R$.
Replacing $R$ by an iterate if necessary we assume that $a = R^\ell(c)$ is a fixed point of $R$.
Denote $g$ a local inverse of $R$ at $a$.

As $c$ is not periodic under $R$, the sets $R^{-n}(c)$ for $n \ge 1$ are pairwise disjoint.
It follows that for large $n \ge 1$ and every $w' \in R^{-n}(c)$, all preimages of $w'$ are unramified.
As $c$ is unramified, there is such $w'$ that is in addition an unramified preimage of $c$.

By the eventually onto property, there is a point $w$ in the domain of definition of $g$ and an integer $m$ such that $R^m(w) = w'$~\cite{Hs} (note that $w'$ is not an exceptional point, as $R^n(w') = a$ is a repelling fixed point.)
By the choice of $w'$, the point $w$ is an unramified preimage of $c$.
So the backward orbit $\{ g^i(w) \}_{i \ge 0}$ and the forward orbit of $w$ form a critical homoclinic orbit.

As $(R, c)$ is defined over $K$ and $w \in R^{-(m + n)}(c)$, it follows that $w$ and its forward orbit belong to a finite extension $L$ of $K$.
As $a = R^\ell(c) \in K$ it follows that $g$ is defined over $K$, and hence over $L$.
So the backward orbit $\{ g^i(w) \}_{i \ge 0}$ is contained $L$.
$\hfill \square$

\subsection{One parameter families of rational maps.}
A {\it one parameter family of rational maps} defined over a ball $B$ of $\C_p$, is a function
$$
R : B \times \pp \to \pp
$$
such that for each parameter $t \in B$ we obtain a rational map $R_t : \pp \to \pp$ {\it whose degree is independent of} $t$ and whose coefficients are holomorphic functions defined on $B$.
We say that a one parameter family of rational maps $R$ {\it is defined over} a finite extension $K$ of $\Q_p$, if the coefficients of $R_t$ are defined over $K$, as holomorphic functions.

\

The following are natural operations to obtain new families from a given one.
\begin{enumerate}
\item[1.]{\it Coordinate change.}
Replace $R_t$ by $\varphi_t \circ R_t \circ \varphi_t^{-1}$, where $\varphi : B \times \pp \to \pp$ is a one parameter family of M\"oebious transformations (rational maps of degree~1).
\item[2.]{\it Iteration.}
For a given $n \ge 1$, replace $R_t$ by $R_t^n$.
\item[3.] {\it Base change.}
Replace $R_t$ by $R_{\Phi(t)}$, where $\Phi : B \to B$ is a non-constant holomorphic function.
\end{enumerate}
Note that these operations preserve families defined over a finite extension $K$ of $\Q_p$, when $\varphi$ and $\Phi$ are defined over $K$.
\subsection*{Marked critical points.}
A {\it one parameter family of rational maps with marked critical point} $(R, c)$ is a one parameter family of rational maps $R : B \times \pp \to \pp$, together with a holomorphic function $c : B \to \pp$ such that $c_t$ is a critical point of $R_t$, for every $t \in B$.
We say that $(R, c)$ is defined over a finite extension $K$ of $\Q_p$, when both $R$ and $c$ are defined over $K$.

\subsection{}
The following lemma shows that a one parameter family of rational maps is locally a one parameter family of holomorphic functions.
So the considerations of Section~\ref{holomorphic functions} apply to families of rational maps.

\

\lemma
{\it Let $R$ be a one parameter family of rational maps defined over a ball $B$.
Then for every $(t_0, z_0) \in B \times \C_p$ such that $z_0$ is not a pole of $R_{t_0}$, the restriction of $R$ to a sufficiently small bidisk containing $(t_0, z_0)$ is a one parameter family of holomorphic functions.}
\proof
After a coordinate change assume $z_0 = R_0(z_0) = 0$.
As by assumption the degree of $R_t$ does not depend on $t$, we can write $R_t$ as the quotient of polynomials $R_t = P_t / Q_t$, in such a way that the coefficients of $P_t$ and of $Q_t$ are holomorphic functions of $t$ and such that for each $t \in B$ the polynomials $P_t$ and $Q_t$ do not have common factors.
Write $Q_t(z) = c_0(t) + c_1(t) z + \ldots + c_d(t) z^d$.
As by assumption $P_0$ and $Q_0$ do not have common factors and $R_0(0) = 0$, we have $P_0(0) = 0$ and $c_0(0) \neq 0$.
Shrinking $B$ if necessary, assume that the function $t \mapsto c_0(t)^{-1}$ is well defined and holomorphic on $B$.
Then, replacing $P$ by $c_0 P$ and $Q$ by $c_0^{-1}Q$ if necessary, we assume $c_0 \equiv 1$.

Let $h$ be the holomorphic function $1 - Q$ and let $\underline{D}$ be a sufficiently small polydisk so that $\sup_{\underline{D}} |h| < 1$.
Then, on $\underline{D}$ we have that $R = P(1 + h + h^2 + \ldots)$ is a holomorphic function.
$\hfill \square$
\subsection{Holomorphic dependence.}\label{holomorphic dependence}
Let $R$ be a one parameter family of holomorphic functions defined over a ball $B$.
Let $t_0 \in B$ and $a_0 \in \C_p$ be such that $a_0$ is a repelling fixed point of $R_{t_0}$.
Then there is a holomorphic function $a$, defined on a ball $B'$ containing $t_0$, such that for every $t \in B'$ the point $a_t$ is a repelling fixed point of $R_t$.
The function $\lambda_t = R_t'(a_t)$ is holomorphic.

Let $g$ be a one parameter family of holomorphic functions defined on a bidisk containing $(t_0, a_0)$, such that for every $t$ the holomorphic function $g_t$ is a local inverse of $R_t$ at $a_t$ (Lemma~\ref{one parameter inverse branches}).
Then we can show in a similar way as in Lemma~\ref{periodic points}, that the sequence of holomorphic functions $\{ \lambda^i (g^i - a) \}_{i \ge 0}$ converges uniformly to a one parameter family of holomorphic functions $\zeta$.

When $R$ is defined over a finite extension $K$ of $\Q_p$ and $t_0, a_0 \in K$, then $a$, $\lambda$, $g$ and $\zeta$ are all defined over $K$.
\section{Misiurewicz bifurcations.}\label{Misiurewicz bifurcations}
Let $(R, c)$ be a one parameter family of rational maps with marked critical point.
We say that $(R, c)$ has a {\it Misiurewicz bifurcation at $t = t_0$}, if the following conditions are satisfied.
\begin{enumerate}
\item[(M1)]
There is an integer $\ell \ge 1$ such that $R_{t_0}^\ell(c_{t_0})$ is a repelling fixed point of
$R_{t_0}$.
\item[(M2)]
$c_{t_0}$ is unramified for $R_{t_0}$.
\item[(M3)]
$R_t^\ell(c_t)$ is not a fixed point of $R_t$, for some parameter $t$.
\item[(M4)]
$\deg_{R_t^\ell}(c_t) = d$ for every $t$ close to $t_0$.
\end{enumerate}

For $t$ close to $t_0$, denote by $a_t$ the fixed point of $R_t$ that is the continuation of the repelling fixed point $R_{t_0}^\ell(c_{t_0})$ of $R_{t_0}$.
We say that the Misiurewicz bifurcation of $(R, c)$ at $t = t_0$ is {\it transversal}, if the function $t \mapsto R_t^{\ell}(c_t) - a_t$ has a non-zero derivative at $t = t_0$, in a coordinate such that $a_{t_0} \neq \infty$.
\subsection{Admissible homoclinic orbits.}
Let $(R, c)$ be a one parameter family of rational maps with marked critical point, having a Misiurewicz bifurcation at $t = 0$.
Denote by $\ell$ the least integer for which $a_0 = R_0^\ell(c_0)$ is a repelling fixed point of $R_0$.
Assume $a_0 \neq \infty$ and set $\lambda_0 = R_0'(a_0)$.

For a homoclinic orbit $\{ z_i \}_{i \in \Z}$ the limit
$$
\zeta = \lim_{i \to \infty} \lambda_0^i(z_{- i} - a_0)
$$
exists and it is non-zero, cf. Lemma~\ref{periodic points}.
We call $\zeta$ {\it the asymptotic position} of the homoclinic orbit $\{ z_i \}_{i \in \Z}$.
Note that a coordinate change $\varphi$ affects $\zeta$ by a factor of $\varphi_0'(a_0)$.

Set $d = \deg_{R_0^\ell}(c_0)$ and let $\rho \in \C_p^*$ be defined by
$$
R_0^\ell (z) = a_0 + \rho (z - c_0)^d + \ldots
$$
For small $t$ denote by $a_t$ the fixed point of $R_t$ that is the continuation of $a_0$.
As $(R, c)$ has a Misiurewicz bifurcation at $t = 0$, it follows that the function $t \mapsto R_t^\ell(c_t) - a_t$ vanishes at $0$ and that it is non-constant.
So, there is a positive integer $d' \ge 1$ and $\rho' \in \C_p^*$ such that
$$
R_t^\ell(c_t) - a_t = \rho' t^{d'} + \ldots
$$

When $(R, c)$ is defined over a finite extension $K$ of $\Q_p$, then $a_0 = R_0^\ell(c_0), \lambda_0, \rho$ and $\rho'$ all belong to $K$.
If moreover the homoclinic orbit $\{ z_i \}_{i \in \Z}$ is contained in $K$, then $\zeta$ belongs to $K^*$.

\

\definition
{\it
Let $(R, c), \lambda_0, \ldots$ be as above and suppose that $(R, c)$ is defined over a finite extension $K$ of $\Q_p$.
A $K$-{\bf admissible homoclinic orbit} is a critical homoclinic orbit $\{ z_j \}_{j \in \Z}$ contained in $K$, for which there is an integer $m_0$ such that
$$
- \lambda_0^{m_0} \zeta / \rho \in K^d
\ \ \mbox{ and } \ \
\lambda_0^{m_0} \zeta / \rho' \in K^{d'}.
$$}

\

\remark \

\begin{enumerate}
\item[1.]
The $K$-admissibility condition is invariant under coordinate changes and base
changes defined over $K$.
In fact, a coordinate change $\varphi$ affects the first element by a
factor of $(\varphi_0'(c_0))^d$ and leaves the second invariant, while a base change
$\Phi$ leaves the first invariant and affects the second by a factor of
$(\Phi'(0))^{d'}$.
\item[2.]
Every critical homoclinic orbit can be made admissible by a finite extension of the base field.
\item[3.]
When the Misiurewicz bifurcation $(R, c)$ is transversal, that is when $d' = 1$, the second condition of $K$-admissibility is automatically
satisfied.
\end{enumerate}

\subsection{Misiurewicz cascades.}\label{Misiurewicz cascade} \
In this section we reduce the proof of Theorem~A to the following proposition.
The proof of this proposition is deferred to the next section.

\

{\par\noindent {\bf Proposition~A.}}
{\it Let $(R, c)$ be a one parameter family of rational maps with marked critical point, having a Misiurewicz bifurcation at $t = 0$.
Suppose that $(R, c)$ is defined over a finite extension $K$ of $\Q_p$ and that there exists a $K$-admissible homoclinic orbit.

Then there is a sequence of parameters $\{ t_r \}$ in $K$ converging to $0$, such that for every $r$, $(R, c)$ has a transversal Misiurewicz bifurcation at $t = t_r$ that admits a $K$-admissible homoclinic orbit.
Moreover, if we denote by $\ell_r \ge 1$ the least integer such that $R^{\ell_r}(c_{t_r})$ is a repelling fixed point of $R_{t_r}$, then}
$$
\lim_{r \to \infty} \min_{1 \le i \le \ell_r} \Delta(c_{t_r}, R^i_{t_r}(c_{t_r})) = 0.
$$

\

The following corollary is immediate consequence of the proposition.

\

\noindent
{\bf Corollary 1.} (Transversal Misiurewicz Bifurcations){\bf .} 
{\it Let $(R, c)$ be a one parameter family of rational maps defined over a finite extension $K$ of $\Q_p$.
Among the parameters in $K$ for which $(R, c)$ has a Misiurewicz bifurcation admitting a $K$-admissible homoclinic orbit, those having a transversal Misiurewicz bifurcation are dense.}

\

\subsection*{Proof of Theorem~A}
Let $(R, c)$ and $K$ be as in the proposition.
Then there is a critical homoclinic orbit contained in a finite extension $L$ of $K$ (Proposition~\ref{critical homoclinic}).
Taking $L$ larger if necessary, this critical homoclinic orbit is $L$-admissible.
Then Theorem~A is an immediate consequence of the following corollary of the proposition.

\

\noindent
{\bf Corollary 2.}
{\it Let $(R, c)$ be a Misiurewicz bifurcation at $t = 0$, defined over $K$.
If there exists a $K$-admissible homoclinic orbit, then there are arbitrarily small parameters $t$ in $K$ for which the critical point $c_t$ of $R_t$ is recurrent, but not periodic.}
\proof
Let $m_0$ be the least positive integer such that $a_0 = R_0^{m_0}(c_0)$ is a repelling fixed point of $R_0$. 
For small $t$ denote by $a_t$ the repelling fixed point of $R_t$ that is the continuation of $a_0$.

Set $s^0 = 0$.
From the proposition we can define by induction a sequence of parameters $\{ s_n \}_{n \ge 1}$ in $K$, converging to some parameter $t_0 \in K$, such that for every $n \ge 1$ the following properties hold.
\begin{enumerate}
\item[1.]
 $(R, c)$ has a Misiurewicz bifurcation at $t = s_n$ that admits a $K$-admissible homoclinic orbit.
\item[2.]
If we denote $m_n$ the least integer such that $R_{s_n}^{m_n}(c_{s_n}) = a_{s_n}$, then
$$
\min_{1 \le i \le m_n} \Delta(R_{s_n}^i(c_{s_n}), c_{s_n}) \le 2^{-n}.
$$
\item[3.]
$s_n$ is close enough to $s_{n - 1}$ so that for every $1 \le i \le m_{n - 1}$ we have
$$
\Delta(R_{s_n}^i(c_{s_n}), R_{s_{n - 1}}^i(c_{s_{n - 1}})) \le 2^{-n}.
$$
\end{enumerate}

Property~2 implies that for every $n \ge 1$ we have
$$
\min_{1 \le i \le m_n} \Delta(R_{t_0}^i(c_{t_0}), c_{t_0}) \le 2^{-n},
$$
so that the critical point $c_{t_0}$ of $R_{t_0}$ is recurrent.
Property~3 implies that for every $n \ge 1$ we have
$$
\Delta(R_{t_0}^{\ell_n}(c_{t_0}), a_{t_0}) \le 2^{-n},
$$
so that the forward orbit of $c_{t_0}$ under $R_{t_0}$ accumulates on $a_{t_0}$.
As $a_{t_0}$ is a repelling fixed point of $R_{t_0}$, it follows that the critical point $c_{t_0}$ is not periodic.
$\hfill \square$

\subsection{Recurrent critical points in $\Q_p$.}
Given a prime number $p$ and an integer $d > 1$, consider the one parameter family of polynomials
$$
P_t(z) = - p^{-d} z(z - 1)^d + t.
$$
As for every parameter $t$ the local degree of $P_t$ at $1$ is equal to $d > 1$, it follows that $(P, 1)$ is a one parameter family of polynomials with marked critical point, that is defined over $\Q_p$.
We will verify that this family has a Misiurewicz bifurcation at $t = 0$ and that there is a $\Q_p$-admissible homoclinic orbit.
Corollary~2 then implies that there are arbitrarily small parameters $t$ in $\Q_p$ for which the critical point $1$ of $P_t$ is recurrent, but not periodic.

Let us now verify that $(P, 1)$ has a Misiurewicz bifurcation at $t = 0$.
\begin{enumerate}
\item[(M1)]
$P_0(1) = 0$ is a fixed point of $P_0$ whose multiplier $\lambda_0 = -(-1)^d p^{-d}$ satisfies $|\lambda_0| = |p|^{-d} > 1$.
So property (M1) is satisfied with $\ell = 1$.
\item[(M2)]
Note that $P_0$ is injective on $\fm_p = \{ |z| < 1 \}$ and $P_0(\fm_p) = \{ |z| < |p|^{-d} \}$.
So, if we denote by $g$ the inverse of $P_0|_{\fm_p}$, then $\{ g^i(1) \}_{i \ge 1}$ is an infinite sequence of unramified preimages of $1$.
This proves that~1 is unramified for $P_0$.
\item[(M3)]
As $P_t(0) = t$, for $t \neq 0$ the point $P_t(1) = 0$ is not fixed by $P_t$.
\item[(M4)]
As $\ell = 1$, for every $t$ we have $\deg_{P_t^\ell}(1) = d$.
\end{enumerate}
Next observe that this Misiurewicz bifurcation is transversal.
In fact, for small $t$ denote by $a_t$ the repelling fixed point of $P_t$ that is the continuation of $a_0 = 0$.
From the equation $a_t = - p^{-d}a_t(a_t - 1)^d + t$ we deduce that $|a_t| = |p^d t|$ and that the derivative of $t \mapsto P_t(1) - a_t$ at $t = 0$ is different from zero.

We will show now that the homoclinic orbit formed by the backward orbit $\{ g^i(1) \}_{i \ge 1}$ and the forward orbit of~1 under $P_0$ is $\Q_p$-admissible.
As our Misiurewicz bifurcation is transversal, we only need to verify the first condition of $\Q_p$-admissibility.

Since $P_0$ and $g$ are defined over $\Q_p$, this homoclinic orbit is contained in $\Q_p$.
Denote by
$$
\zeta = \lim_{i \to \infty} \lambda_0^i g^i(1),
$$
its asymptotic position.
As
$$
P_0(z) = -p^{-d}(z - 1)^d - p^{-d} (z - 1)^{d + 1},
$$
we have $\rho = -p^{-d}$.
We will show that $\zeta \in (\Q_p)^d$, which implies that $- \zeta / \rho = p^d \zeta \in (\Q_p)^d$ and that the first condition of $\Q_p$-admissibility is satisfied with $m_0 = 0$.
Set $\zeta_0 = 1$ and for $i \ge 1$ set $\zeta_i = \lambda_0^i g^i(1)$, so that $\zeta = \lim_{i \to \infty} \zeta_i$.
From the identity $P_0(g(z)) = z$ we have $\lambda_0 g(z)/z = (1 - g(z))^{-d} \in (\Q_p^*)^d$, whenever $z \in \Q_p$ satisfies $|z| \in (0, |p|^{-d})$.
It follows that for $i \ge 0$ we have $\zeta_{i + 1}/\zeta_i \in (\Q_p)^d$ and that $\zeta_{i + 1} / \zeta_i \to 1$ as $i \to \infty$.
Therefore
$$
\zeta = \prod_{i \ge 0} \left( \zeta_{i + 1}/\zeta_i \right) \in \Q_p^d.
$$

\section{Proof of Proposition~A.}
Denote by $\ell \ge 1$ the least integer for which $a_0 = R_0^\ell(c_0)$ is a repelling fixed point of $R_0$.
For small $t$ denote by $a_t$ the repelling fixed point of $R_t$ that is the continuation of $a_0$.
Since every non-constant rational map has at least one non-repelling fixed point~\cite{Be1}, after a coordinate change we assume that $\infty$ is a non-repelling fixed point of $R_0$.
So for every small~$t$ we have $a_t \neq \infty$ and $c_t \neq \infty$.

Denote by $\lambda_t = (R_t^\ell)'(a_t)$ the multiplier of $a_t$ and denote by $g_t$ the local inverse of $R_t$ at $a_t$.
We assume that the $g_t$ are all defined on a ball of $\C_p$ containing $a_0$, that does not depend on $t$ (Section~\ref{one parameter inverse branches}).

\noindent
{\bf A.}
By hypothesis there is a $K$-admissible homoclinic orbit $\{ z_i \}_{i \in \Z}$ of $(R_0, c_0)$; denote by $\zeta$ its asymptotic position.
Let $\rho, \rho' \in K^*$, $d = \deg_{R_t^\ell}(c_t) > 1$ and $d' \ge 1$ be as in the admissibility condition.
Then there are $m_0 \in \Z$ and $\xi, \xi' \in K^*$ such that
$$
\rho \xi^d = - \lambda_0^{m_0} \zeta
\ \mbox{ and } \
\rho' (\xi')^{d'} = \lambda_0^{m_0} \zeta.
$$

For a given integer $r$ set $m = r dd' - m_0$,
$$
B' = \{ z \in \C_p \mid |z - \lambda_0^{-rd} \xi'| \le |p^2d| \cdot |\lambda_0^{-rd} \xi'| \}
$$
and for $t \in B'$,
$$
B_t = \{ z \in \C_p \mid |z - \lambda_0^{-rd'} \xi - c_t| \le |p^2d'| \cdot |\lambda_0^{-rd'} \xi| \}.
$$
Lemma~\ref{near critical inverse branches} applied to $\delta = |p^2d|$ and to $t \mapsto R_t^\ell(c_t) - a_t$, implies that when $r$ is large enough the function $t \mapsto R_t^\ell(c_t) - a_t$ maps $B'$ bijectively onto
$$
\{ w \in \C_p \mid |w - \lambda_0^{-m} \zeta| \le |p^2dd'| \cdot |\lambda_0^{- m} \zeta| \}.
$$
Moreover, Lemma~\ref{near critical inverse branches} applied to $\delta = |p^2d'|$ and to the family $z \mapsto R_t^\ell(z + c_t) - R_t^\ell(c_t)$, implies that when $r$ is large enough for every $t \in B'$ the function $R_t^\ell$ maps $B_t$ bijectively onto
$$
\{ w \in \C_p \mid |w - R_t^\ell(c_t) + \lambda_0^{- m} \zeta| \le |p^2dd'| \cdot |\lambda_0^{- m} \zeta| \}
$$
$$
=
\{ w \in \C_p \mid |w - a_t| \le |p^2dd'| \cdot |\lambda_0^{- m} \zeta| \}.
$$

\noindent
{\bf B.}
As the homoclinic orbit $\{ z_i \}_{i \in \Z}$ is critical, we have $z_0 = c_0$ and for every $j > 0$ the point $z_{-j}$ is an unramified preimage of $c_0$.
Moreover, for some $n > 0$ the point $z_{- n}$ belongs to the domain of $g_0$ and for $j \ge 0$ we have $z_{- (n + j)} = g_0^j(z_{- n})$.
Since $\infty$ is a non-repelling fixed point of $R_0$, the homoclinic orbit $\{ z_i \}_{i \in \Z}$ is contained in $\C_p$.

By Lemma~\ref{one parameter inverse branches} it follows that, for every $0 \le i \le n$ and every sufficiently small ball $\cD$ containing $c_0$, there is a ball $\cB$ containing $t_0$ and a one parameter family of holomorphic functions $h^i$ defined on $\cB \times \cD$, such that for every $t \in \cB$ the holomorphic function $h_t^i$ is a local inverse of $R_t$ such that $h^i_0(c_0) = z_{- i}$.
Reducing $\cD$ if necessary, we assume that for every $0 \le i \le n$ and $t \in \cB$ the set $h^i_t(\cD)$ is contained in the domain of $g_t$.
Then for $j > 0$ the function $h_t^{n + j} = g_0^j \circ h_t^n$ is a holomorphic function defined on $\cB \times \cD$.
So for every $j > 1$ the function $h^j$ is a holomorphic function defined on $\cD \times \cB$ and for every $t \in \cB$ the holomorphic function $h^j_t$ is a local inverse of $R_t^j$ such that $h_0^j(c_0) = z_{ -j}$.

It follows that for every small $t$ the critical point $c_t$ of $R_t$ is unramified, as $\{ h_t^i(c_t) \}_{i \ge 0}$ is an infinite sequence of unramified preimages of $c_t$.

\noindent
{\bf C.}
The proof is organized as follows.
In part~D we define for every $t \in B' = B'(r)$, an orbit $\{ z_i(t) \}_{\ell \le i \le \ell_r}$ ending at $z_{\ell_r}(t) = a_t$.
In part~E we show that when $r$ is sufficiently large there is a parameter $t_r \in B'$ in $K$ such that $R_{t_r}^\ell(c_{t_r}) = z_\ell(t_r)$, so that
$$
R_{t_r}^{\ell_r}(c_{t_r}) = z_{\ell_r}(t_r) = a_{t_r}
$$
and so that $(R, c)$ has a Misiurewicz bifurcation at $t = t_r$.
In part~F we show this Misiurewicz bifurcation is transversal.
In part~G we show that the critical homoclinic orbit (for $t = t_r$) determined by the infinite backward orbit $\{ h_{t_r}^i(c_{t_r}) \}_{i \ge 0}$, is $K$-admissible.
The last assertion of the proposition is easy to verify, this is done in part~E.

\noindent
{\bf D.}
Let $\hr > r$ be such that $|\lambda_0|^{-(\hr - r)dd'} < |p^2dd'|$ and set $\hm = \hr dd' - m_0$.
Given an positive integer $k$ set $\ell_r = (k - 1)(\hm + \ell) + m + 2 \ell$. 

We assume $r$ large enough so that for every small parameter $t$ we have $B_t \subset \cD$.
For $t \in B'$ define $w_t^1, \ldots, w_t^k \in B_t$ inductively as follows.
Let $w_t^1$ be the unique inverse image of $a_t$ by $R_t^\ell$ in $B_t$.
Let $1 < j \le k$ be given and assume that $w^{j - 1}_t$ is already defined.
When $r$ is sufficiently large we have
$$
|h_t^{\hm}(w^{j - 1}_t) - a_t|
=
|\lambda_0|^{-\hm} |\zeta|
<
|p^2dd'| \cdot |\lambda_0|^{- m}|\zeta|.
$$
So $h_t^{\hm}(w_t^{j - 1}) \in R_t^\ell(B_t)$ has a unique preimage by $R_t^\ell$ in $B_t$.
We let $w_t^j$ be this preimage.

Set $z_\ell(t) = h_t^m(w_t^k)$ and for $\ell < i \le \ell_r$ set $z_i(t) = R_t^{i - \ell}(z_1(t))$.
Then, for $0 \le j \le k - 1$ we have
$$
z_{j(\hm + \ell) + m + \ell}(t) = w_t^{k - j},
$$
and therefore $z_{\ell_r}(t) = a_t$.

It follows from Lemma~\ref{near critical inverse branches} and from the fact that the functions $h^{\nu}$ are holomorphic and defined over $K$, that the functions $w^j$ and $z_i$ are holomorphic functions defined over $K$.

\noindent
{\bf E.}
Recall that $m = rdd' - m_0$ depends on $r$.
Assume $r$ large enough so that for every $w$ close to $c_t$ and every small~$t$ we have
$$
\left| \zeta - \lambda_0^m(h_t^m(w) - a_t) \right|
\le
|p^3 dd'| \cdot |\zeta|
$$
(Section~\ref{holomorphic dependence}).
Since for every $t \in B' \subset \{ |z| \le |\lambda_0|^{- rd}|\xi'| \}$ we have $B_t \subset \{ |z| \le |\lambda_0|^{- rd'}|\xi| \}$, it follows that for large $r$ the function
$$
g(t)
=
\lambda_0^{- m} \zeta - (h_t^m(w_t^k) - a_t),
$$
satisfies $|g| \le |p^3 dd'| \cdot |\lambda_0|^{- m}|\zeta|$ on $B'$.

On the other hand, by part~A the function $f(t) = R_t^\ell(c_t) - a_t - \lambda_0^{- m} \zeta$ maps the ball $B'$ univalently onto the ball
$$
\{ w \in \C_p \mid |w| \le |pdd'| \cdot |\lambda_0|^{- m} |\zeta| \}.
$$
So $\sup_{B'} |g| < \diam(f(B'))$ and Lemma~\ref{Rouche} implies then that there is parameter $t_r \in B'$ such that $f(t_r) = g(t_r)$, so that $R_{t_r}^\ell(c_{t_r}) = h_{t_r}^m(w_{t_r}^k) = z_\ell(t_r)$.
Since $(R, c)$ and the functions $a$ and $h^m(w^k)$ are defined over $K$ and $\lambda_0^{-m} \zeta \in K$, it follows that $f$ and $g$ are defined over $K$ and $t_r \in K$ (cf. Lemma~\ref{Rouche}).

Let us verify now that $(R, c)$ has a Misiurewicz bifurcation at $t = t_r$.
By definition of $t_r$ we have $R_{t_r}^{\ell_r}(c_{t_r}) = z_{\ell_r}(t_r) = a_{t_r}$.
As remarked in part~B, for every $t \in B'$ the critical point $c_t$ is unramified for $R_t$.
By construction, for every $t \in B'$ none of the points $z_i(t)$, for $\ell \le i \le \ell_r$, is a critical point of $R_t$, so $\deg_{R^{\ell_r}}(c_t) = \deg_{R^{\ell}} (c_t) = d$.
It remains to verify condition (M3). 
Suppose by contradiction that for every $t \in B'$ we have $R_t^{\ell_r}(c_t) = a_t$.
This implies that the functions $R_t^\ell(c_t)$ and $z_\ell(t)$, and hence $f$ and $g$, coincide on $B'$.
But we saw above that $g(B')$ is strictly smaller than $f(B')$.

To verify the last assertion of the proposition, observe that $R_{t_r}^{m + \ell}(c_{t_r}) = w_{t_r}^k \in B_{t_r}$, so that
$$
\Delta(c_{t_r}, R_{t_r}^{m + \ell}(c_{t_r})) \le |\lambda_0|^{- rd}|\xi| \to 0
\mbox{ as }
r \to \infty.
$$

\noindent
{\bf F.}
We prove now that for large $r$ the Misiurewicz bifurcation $(R, c)$ at $t = t_r$ is transversal.
As $\infty$ is a non-repelling fixed point of $R_0$,  it follows that for small $t$ the forward orbit of $c$ under $R_{t_r}$ does not contain $\infty$.
So $\partial_t R_t$ is well defined at each of these points.

After an affine change of coordinates we assume that $a = a_t$ and $c = c_t$ do not depend on $t$.
So for all $t$ we have $R_t(a) = a$ and $\partial_t R_t (a) = 0$.
Moreover there is a constant $C > 0$ such that for every point $z$ close to $a$ we have $|\partial_t R_t(z)| \le C |z - a|$.
Enlarging $C$ if necessary we assume that for every $w \in \cD$ near $c$ and every positive integer $i$ we have
\begin{equation}\label{transversality 1}
\left| \partial_t (R_t)|_{t = t_r} (h_{t_r}^i(w)) \right| \le C |\lambda_0|^{-i}
\ \mbox{ and } \
|(R_{t_r}^i)'(h_{t_r}^i(w))| \ge C^{-1} |\lambda_0|^i.
\end{equation}

\noindent
$\bullet$
For $0 \le j \le k - 1$ set $N_j = m + \ell + j(\hm + \ell)$, so that $\ell_r = N_{k - 1} + \ell$.
Moreover for every $0 \le j \le k - 1$ we have $R_{t_r}^{N_j}(c_{t_r}) = w_{t_r}^{k - j}$ and for $1 \le j \le k - 1$ we have $R_{t_r}^{N_{j - 1} + \ell}(c_t) = h_{t_r}^{\hm}(w^{k - j})$.

For $0 \le j \le k - 1$ set $D_j = \partial_t \left( R_t^{N_j + \ell}(c_t) \right)|_{t = t_r}$.
Note that $D_{k - 1} = \partial_t (R_t^{\ell_r}(c_t) - a)$, so we need to prove that $D_{k - 1} \neq 0$.

For $0 \le j \le k - 1$ we have
$$
D_j = A_j + B_j + C_j,
$$
where $A_j = (\partial_t R_t^\ell)|_{t = t_r}(w_{t_r}^{k - j})$, and when $j \neq 0$,
$$
B_j = (R_{t_r}^\ell)'(w_{t_r}^{k - j}) \cdot \partial_t (R_t^{\hm})|_{t = t_r}(h_{t_r}^{\hm}(w_{t_r}^{k - j})),
$$
$$
C_j = (R_{t_r}^{\hm + \ell})'(h_{t_r}^{\hm}(w_{t_r}^{k - j})) \cdot \partial_t \left(R_t^{N_{j - 1} + \ell}(c_t) \right)|_{t = t_r}.
$$
When $j = 0$,
$$
B_0 = (R_{t_r}^\ell)'(w_{t_r}^k) \cdot \partial_t (R_t^{m})|_{t = t_r}(h_{t_r}^{m}(w_{t_r}^k)),
$$
$$
C_0 = (R_{t_r}^{m + \ell})'(h_{t_r}^m(w_{t_r}^{k})) \cdot \partial_t \left(R_t^{\ell}(c_t) \right)|_{t = t_r}.
$$

We will show by induction in $j$ that $|D_j| = |C_j|$.
By induction in $j$ this implies that
$$
|D_j|
= 
\left| (R_{t_r}^{N_j})'(R_{t_r}^\ell(c_{t_r})) \cdot \partial_t \left(R_t^{\ell}(c_t) \right)|_{t = t_r} \right|.
$$
But this last quantity is non-zero, since the first factor is non-zero by construction and $\left| \partial_t \left(R_t^{\ell}(c_t) \right)|_{t = t_r} \right| = |\rho'd'(t_r)^{d' - 1}| \neq 0$.
So this claim implies that $D_{k - 1} \neq 0$, as wanted.

\noindent
$\bullet$
Observe that for $t$ close to $0$ and $w$ close to $c$ we have
$$
R_t^\ell(w) = a_0(t) + a_d(t)(w - c)^d + \ldots,
$$
with $a_0(t) = a + \rho't^{d'} + \ldots$ and $a_d(0) = \rho$.
Therefore, enlarging the constant $C$ if necessary, we have
$$
\left |\partial_t (R_t^\ell)|_{t = t_r} (w) \right|
\le
C \max\{ |t_r|^{d' - 1}, |w - c|^d \}.
$$
As $|t_r| = c_0' |\lambda_0|^{-m/d'}$ and $|w_{t_r}^{k - j} - c| = c_0|\lambda_0|^{-m/d}$ for constants $c_0, c_0' > 0$ independent of $j$ and $r$, it follows that when $r$ is large enough we have
$$
|A_j| \le C (c_0')^{d' - 1} |\lambda_0|^{-m(1 - 1/d')}.
$$

\noindent
$\bullet$
For large $r$ we have
$$
|(R_{t_r}^\ell)'(w_{t_r}^{k - j})| = |\rho (w_{t_r}^{k - j} - c)^{d - 1}| = c_1 |\lambda_0|^{-m(1 - 1/d)},
$$
where $c_1 = |\rho|c_0^{d - 1}$.

Fix $0 \le j \le k - 1$ and let $n = \hm$ if $j \neq 0$ and $n = m$ if $j = 0$.
Then
$$
\partial_t \left( R_t^{n} \right)|_{t = t_r} (h_{t_r}^n(w_{t_r}^{k - j}))
=
\sum_{1 \le i \le n}
(R_t^{i - 1})'(h_{t_r}^{i - 1}(w_{t_r}^{k - j}))
\cdot
\partial_t (R_t)|_{t = t_r} (h_{t_r}^{i}(w_{t_r}^{k - j})).
$$
The inequalities ($\ref{transversality 1}$) then imply
$$
\left| \partial_t \left( R_t^{n} \right)|_{t = t_r} (h_{t_r}^n(w_{t_r}^{k - j}))  \right|
\le C^2|\lambda_{0}|^{-1} \le C^2.
$$

So $|B_j| \le C^2c_1 |\lambda_0|^{-m(1 - 1/d)}$.

\noindent
$\bullet$
We have
$$
C_0 = (R_{t_r}^m)'(h_{t_r}^m(w^k)) \cdot (R_{t_r}^\ell)'(w^k) \cdot \partial_t (R_t^\ell(c_t))|_{t = t_r},
$$
so $|C_0| \ge C^{-1}|\lambda_0|^m \cdot c_1 |\lambda_0|^{-m(d - 1)/d} \cdot |\rho'd'| |t_r|^{d' - 1} = c_2 |\lambda_0|^{- m (1 - 1/d - 1/d')}$, where $c_2 = C^{-1} c_1 |\rho'd'|(c_0')^{d' - 1}$.

We assume $r$ large enough so that
\begin{equation}\label{transversality 2}
c_2|\lambda_0|^{- m(1 - 1/d - 1/d')} > \max\{ C(c_0')^{d' - 1} |\lambda_0|^{- m(1 - 1/d') }, C^2c_1 |\lambda_0|^{-m(1 - 1/d)} \}.
\end{equation}
This implies $|C_0| > \max \{ |A_0|, |B_0| \}$ and $|D_0| = |C_0|$.

For $1 \le j \le k - 1$ we have
$$
C_j = (R_{t_r}^{\hm})'(h_{t_r}^{\hm}(w_{t_r}^{k - j})) \cdot (R_{t_r}^\ell)'(w_{t_r}^{k - j}) \cdot D_{j - 1},
$$
so for large $r$ we have $|C_j| \ge C^{-1}|\lambda_0|^{\hm} c_1 |\lambda_0|^{-m(1 - 1/d)} |C_{j - 1}| \ge |C_{j - 1}|$.
It follows by induction in $j$ that $|C_j| \ge |C_0|$.
Then inequality ($\ref{transversality 2}$) implies that $|C_j| > \max \{ |A_j|, |B_j| \}$ and $|D_j| = |C_j|$.

\noindent
{\bf G.}
Now we will prove that, if $k = d$ and $r$ is large enough, then the critical homoclinic orbit, for $t = t_r$, determined by the infinite backward orbit $\{ h_{t_r}^i(c_{t_r}) \}_{i \ge 0}$ is $K$-admissible.
Denote by $\zeta_r$ the asymptotic position of this homoclinic orbit.

By part~F the Misiurewicz bifurcation $(R, c)$ at $t = t_r$ is transversal, so we only need to verify the first condition of $K$-admissibility.
We fix the parameter $t = t_r$ and omit all the indexes corresponding to parameters.
We will use several times that $(K^*)^d$ is open in $K^*$, so that any element in $K$ close to~1 belongs to $K^d$.

Define $\rho_r, \trho_r \in K^*$ by
$$
R^{\ell_r}(z) = a + \rho_r(z - c)^d + \ldots
\mbox{ and }
R^\ell(z) = R^\ell(c) + \trho_r(z - c)^d + \ldots
$$
We have $\rho_r = (R^{\ell_r - \ell})'(R^\ell(c)) \cdot \trho_r$.
As $\lambda = \lambda_{t_r} \to \lambda_0$, $\zeta_r \to \zeta$ and $\trho_r \to \rho$ as $r \to \infty$, it is enough to prove that
$$
(R^{\ell_r - \ell})'(R^\ell(c_t)) = \rho_r/\trho_r \in K^d.
$$
By construction $R^\ell(c) = h^m(w^k)$.
Moreover, recall that
$$
\ell_r - \ell = (k - 1)(\hm + \ell) + m + \ell
\ \mbox{ and } \
\hm = \hr dd' - m_0 = m + (\hr - r) dd'.
$$
Therefore,
\begin{eqnarray*}
(R^{\ell_r - \ell})'(R^\ell(c))
& = &
(R^{\ell_r - \ell})'(h^m(w^k)) \\
& = & 
(R^{k(\hm + \ell)})'(h^{\hm}(w^k)) \cdot \left( (R^{(\hr - r)dd'})'(h^{\hm}(w^k)) \right)^{-1}.
\end{eqnarray*}
We will verify that each of the two factors of the last term belong to $K^d$, thus completing the proof of the proposition.
Observe that, as $t_r \in K$, we have that $c, a, \lambda, w^j, h^i(w^j), \ldots$ all belong to $K$.

\noindent
$\bullet$
For the second factor, observe that when $r$, and hence $m$, is large enough, for every $m < i \le \hm$ the point $h^i(w^k)$ is close enough to the fixed point $a$, so that $R'(h^i(w^k))/\lambda \in K^d$.
Therefore
$$
(R^{(\hr - r)dd'})'(h^{\hm}(w^k))
=
\lambda^{(\hr - r)dd'} \prod_{m < i \le \hm} (R'(h^i(w^k))/\lambda)
$$
belongs to $K^d$, as it is the product of elements in $K^d$.

\noindent
$\bullet$
For the first factor, observe that
$$
(R^{k(\hm + \ell)})'(h^{\hm}(w^k))
=
\prod_{1 \le j \le k} (R^{\hm + \ell})'(h^{\hm}(w^j)).
$$

Since for every $w \in \cD$ we have $h^i(w) \to a$ as $i \to \infty$, it follows that if $w \in K$ is close enough to $c$, then for every positive integer $s$ we have
$$
R'(h^s(w)) / R'(h^s(c)) \in K^d.
$$
So if $r$ is sufficiently large then for every $1 \le j \le k$ we have
$$
(R^{\hm})'(h^{\hm}(w^j)) / (R^{\hm})'(h^{\hm}(c)) \in K^d.
$$

On the other hand, for $1 \le j \le k - 1$ we have
$$
|h^{\hm}(w^j) - a| / |R^\ell(c) - a| = |\lambda_0|^{-(\hr - r)dd'}.
$$
Since $R^\ell(w^1) = a$ and for $1 \le j \le k - 1$ we have $R^\ell(w^{j + 1}) = h^{\hm}(w^j)$, it follows that, if $\hr$ is large enough compared to $r$, the quotient
$$
(R^\ell)'(w^{j + 1}) / (R^\ell)'(w^1)
$$
is close enough to~1, to belong to $K^d$.
Setting $\eta = (R^{\hm})'(h^{\hm}(c)) \cdot (R^\ell)'(w^1)$, we thus have $(R^{\hm + \ell})'(h^{\hm}(w^j)) \eta^{-1} \in K^d$.
When $k = d$ it follows that
$$
\prod_{1 \le j \le k} (R^{\hm + \ell})'(h^{\hm}(w^j))
=
\eta^d \cdot \prod_{1 \le j \le d} \left( (R^{\hm + \ell})'(h^{\hm}(w^j)) \eta^{-1} \right)
\in K^d.
\hfill \square
$$

\bibliographystyle{plain}

\begin{thebibliography}{FvP}

\bibitem[Be1]{Benowandering}
R. Benedetto.
{\it $p$-adic dynamics and Sullivan's No Wandering Domains theorem}.
Compositio Mathematica {\bf 122} (2000), 281-298. 

\bibitem[Be2]{Be1}
R. Benedetto.
{\it Hyperbolic maps in $p$-adic dynamics.}
Ergod. Th. Dynam. Sys. {\bf 21} (2001), 1-11.

\bibitem[Be3]{components}
R. Benedetto.
{\it Components and periodic points in non-Archimedean dynamics.}
Proc. London Math. Soc. {\bf 83} (2002), 231-256.

\bibitem[Be3]{Bewandering}
R. Benedetto.
{\it Examples of wandering domains in p-adic polynomial dynamics.}
C.R. Acad. Sci. Paris {\bf 335} (2002), 615-620.

\bibitem[BGR]{BGR}
S.~Bosch, U.~G\"untzer, R.~Remmert.
{\it Non-Archimedean analysis. A systematic approach to rigid analytic geometry}.
Springer-Verlag, Berlin, 1984.

\bibitem[FvP]{FvP}
J.~Fresnel, M.~van~der~Put. 
{\it G\'eom\'etrie analytique rigide et applications}. 
Progress in Mathematics, 18.
Birkh\"auser, Boston, Mass., 1981

\bibitem[Hs]{Hs}
Liang-Chung Hsia. 
{\it Closure of periodic points over a non-Archimedean field}. 
J. London Math. Soc. {\bf 62} (2000), 685--700.

\bibitem[Mc1]{Mc1}
C.~McMullen.
{\it Families of rational maps and iterative root-finding algorithms.}
Ann. of Math. {\bf 125} (1987), 467-493.

\bibitem[Mc2]{Mc2}
C.~McMullen.
{\it The Mandelbrot set is universal.}
The Mandelbrot set, theme and variations, 1-17, London Math. Soc. Lecture Note Ser., 274, Cambridge Univ. Press, Cambridge, 2000.
 
\bibitem[R1]{these}
J.~Rivera-Letelier.
{\it Dynamiques des fonctions rationnelles sur des corps locaux.}
Ast\'erisque~{\bf 287} (2003), 147-230.

\bibitem[R2]{Fatou}
J.~Rivera-Letelier.
{\it Sur l'estructure des ensembles de Fatou $p$-adiques.}
Preprint~2002. {\tt www.math.sunysb.edu/$\sim$rivera}.

\bibitem[Su]{Su}
D.~Sullivan.
{\it Quasiconformal homeomorphisms and dynamics. I. Solution of the Fatou-Julia problem on wandering domains.}
Ann. of Math. {\bf 122} (1985), 401--418. 
 
\end{thebibliography}

\end{document}